\documentclass[11pt,a4paper]{article}

\usepackage{graphicx}        
\usepackage{amsmath}
\usepackage{enumitem}
\usepackage{amsfonts}
\usepackage{color}
\usepackage[ruled,vlined]{algorithm2e}

\usepackage[a4paper,left=3.cm,right=3.cm,top=2.cm,bottom=3cm]{geometry}

\DeclareMathOperator{\argmin}{argmin}
\DeclareMathOperator{\proj}{Proj}
\newcommand{\uargmin}[1]{\underset{#1}{\argmin}\;}
\newtheorem{prop}{Proposition}
\newtheorem{theorem}{Theorem}

\newtheorem{proper}{Property}
\newtheorem{lemma}{Lemma}

\newcommand{\R}{\mathbb{R}}

\renewcommand{\H}{\mathcal{H}}

\newtheorem{rem}{Remark}

\newcommand{\prox}{\operatorname{prox}}

\newcommand{\norm}[1]{|\!| #1 |\!|}
\newcommand{\normu}[1]{\norm{#1}_{1}}

\newcommand{\normd}[1]{\norm{#1}_{2}}

\def\R{{\mathbb R}}

\newcommand{\hbindex}[1]{#1}  

\begin{document}

\title{Iterative Methods for Computing Eigenvectors of Nonlinear Operators}
\author{Guy Gilboa \\
\small{Technion - Israel Institute of Technology, {guy.gilboa@ee.technion.ac.il}}}
\date{}
%
%
\maketitle
%
\abstract{In this chapter we are examining several iterative methods for solving nonlinear eigenvalue problems. These arise in variational image-processing, graph partition and classification, nonlinear physics and 
39
 more. The canonical eigenproblem we solve is $T(u)=\lambda u$, where $T:\R^n\to \R^n$ is some bounded nonlinear operator. Other variations of eigenvalue problems are also discussed. We present a progression of 5 algorithms, coauthored in recent years by the author and colleagues. Each algorithm attempts to solve a unique problem or to improve the theoretical foundations. The algorithms can be understood as nonlinear PDE's which converge to an eigenfunction in the continuous time domain. This allows a unique view and  understanding of the discrete iterative process. Finally, it is shown how to evaluate numerically the results, along with some   examples and insights related to priors of nonlinear denoisers, both classical algorithms and ones based on deep networks.}

\section{Introduction and Preliminaries} \label{sec_prelim}
In this section, we outline some basic notations and properties which will be used throughout this chapter.
A main type of functionals we are discussing are \hbindex{one-homogeneous functionals}, used frequently as regularizers in image processing and learning.
\subsection{One homogeneous functionals} 
We consider an absolutely one homogeneous functional $J$ that takes as input a function $u:x\in\Omega\to \mathbb{R}$ defined on a  domain $\Omega\subset\mathbb{R}^2$.
$\Omega$ can either be a discrete domain of size $|\Omega|=N$ or an open convex bounded set with Lipschitz boundary.
$u$ are elements of some Hilbert space $X$ (e.g. $X$ can be $L^2(\Omega)$) embeded with some inner product $\langle . \, , \, \rangle$.
$J: X \to \mathbb{R} \bigcup \{+ \infty\}$
is assumed to be proper, convex and lower semi-continuous (lsc).
Absolutely one-homogeneous functionals satisfy 
\begin{equation}\label{eq:J_1hom}
J(cu)=|c|J(u), \,\,\,\, \forall c\in \mathbb{R},\, \forall u\in X.
\end{equation}
The functional $J$ in finite dimensions can be, for instance, of the general form:
\begin{equation}
\label{J}
J(u)=\sum_{i=1}^N\left(\sum_{j=1}^Nw_{ij}|u_i-u_j|^q\right)^{1/q},
\end{equation}
for $q\geq 1$, with $w_{ij}\geq 0$ (usually symmetric weights are assumed
$w_{ij}=w_{ji}$).
This formulation can be understood as a typical one-homogeneous functional on weighted graphs. In this case $u_i$ is the value of the function $u$ at node $i$ on the  graph and $w_{ij}$ is the weight between node $i$ and node $j$. As grids of any
dimension can be realized by specific graph structures, this formulation applies to standard grids as well. Thus \eqref{J}, with appropriate weights, can be the spatial discrete version of anisotropic \hbindex{total variation} (\hbindex{TV}) ($q=1$), isotropic TV ($q=2$) and anisotropic or isotropic \hbindex{nonlocal TV}. 

We  recall the \hbindex{subgradient} definition for general convex functionals
\begin{equation*}
\begin{split}p\in\partial J(u)&\Leftrightarrow J(v)-J(u)\geq \langle p,v-u\rangle,\,\, \forall v.\\
\end{split}
\end{equation*}
We also note the relation to the convex conjugate $J^*$
$$J(u)=\sup_p \langle u, p\rangle-J^*(p).$$
Below we state some properties of one-homogeneous functionals.
\begin{proper}
A function $J$ defined in \eqref{J} admits:
\begin{itemize}
\item [(a)] If $p\in \partial J(u)$, then $J(u)=\langle p,u\rangle$,
\item [(b)]  If $p\in \partial J(u)$, then $J(v)\geq \langle p,v\rangle$, $\forall v$. 
\end{itemize}
\end{proper}
Notice in particular that from (b) we get that $\partial J(u) \subset \partial J(0)$ $\forall u \in X$.

\begin{proper}
The convex conjugate $J^*$ of a one-homogeneous functional
 is the characteristic function of the convex set $\{\partial J(0)\}$.
Moreover, when $\Omega$ is included in a finite dimensional space, we have (\cite{Burger16}):
\begin{equation}\label{bound_subdiff}
\exists C>0\textrm{ s.t. }\normd{p}\leq C,\, \forall p\in\partial J(0).
\end{equation}
\end{proper}
From 
the equivalence of norms, we have that if 
$u$
is of zero mean, there exists a constant $\kappa>0 $ for which
\begin{equation}\label{poincare}
\normd{u}\leq  \kappa  J(u),\,\,\forall u\,\textrm{ such that } \langle u,\mathbf{1}\rangle=0.
\end{equation}

The nullspace of the functional is defined by
\begin{equation}\label{nullspace}
\mathcal{N}(J)=\left\{ u\in X\, |\, J(u)=0\right\}.
\end{equation}
The properties below are shown in \cite{Burger16}.
\begin{proper}
An absolutely one-homogeneous functional $J$ is a seminorm and its nullspace is a linear subspace. 
\end{proper}

\begin{proper}
\label{prop:p}
If a unit constant function $u=\mathbf{1}$ is in $\mathcal{N}(J)$ then any subgradient $p$ admits 
$$\langle p,\mathbf{1}\rangle=0.$$
\end{proper}




We use $\ell_2$ and $\ell_1$ norms of $u$ defined as $\normd{u}=\sqrt{\langle u,u\rangle}$ and $\normu{u}={\langle u,\textrm{sign}(u)\rangle}$.

\section{Eigenvectors of nonlinear operators} 
\label{sec_nlef}
We give here a brief introduction to the broad topic of eigenvectors of nonlinear operators. More details are provided in relation to the variational setting.
We would like to extend the linear eigenvalue problem 
$$ Lu = \lambda u,$$
given a matrix $L$, to a generalized problem, given a bounded nonlinear operator $T:X \to X$.
Replacing $L$ by $T$ we get the \hbindex{nonlinear eigenvalue problem} associated with $T$,
\begin{equation}
\label{eq:ef}
    T(u) = \lambda u,
\end{equation}
where $\lambda \in \mathbb{R}$ is the associated eigenvalue.
In the variational context, given a convex functional $J$, the eigenvalue problem induced by $J$ is 
\begin{equation}
\label{eq:ef_1hom}
    p = \lambda u,   \,\,\,\, p \in \partial J(u).
\end{equation}
As an example, for the Dirichlet energy $J=\frac{1}{2}\|\nabla u\|^2$, the associated eigenvalue problem is a linear one,
$$     -\Delta u = \lambda u,   $$
where $\Delta$ denotes the Laplacian. For appropriate boundary conditions, sines and cosines are solutions to this problem, which are the basis elements of the Fourier transform. For one-homogeneous regularizing functionals, such as total-variation, one obtains different (sharp) eigenfunctions, which can serve for representing signals based on nonlinear \hbindex{spectral} transforms, as shown in \cite{gilboa2013spectral:25,Gilboa_spectv_SIAM_2014,burger2016spectral,gilboa2018book,bungert2019nonlinear}. We would not elaborate on this direction, which is beyond the scope of this chapter.

For absolutely one homogeneous functionals, the eigenvalues are non-negative, since
 $J(u)=\langle\lambda u,u\rangle=\lambda\normd{u}^2$ and $\lambda=\frac{J(u)}{\normd{u}^2}\geq 0$.
An interesting insight on the eigenvalue $\lambda$ shown in \cite{aujol2018theoretical} can be gained by the following proposition. We define $K=\{\partial J(0)\}$ to be the set of possible subgradients for any $u$. Indeed if $p\in\partial J(u)$ then $p\in\partial J(0)$.
We first note that an eigenfunction that admits $\lambda u\in\partial J(u)$ has zero mean from Property \ref{prop:p} above.
Next, we have the following result.
\begin{prop}\label{prop:orth_proj}
For any non constant eigenfunction $u$,  we have $\forall \mu\geq\lambda$,
$$\lambda u=\proj_K(\mu u),$$ where $\proj_K$ is the orthogonal projection onto $K=\{\partial J(0)\}$.
\end{prop}

Eigenfunctions in the form of \eqref{eq:ef_1hom} have analytic solutions, when used as initial conditions in gradient flows. Let
a gradient flow be defined by,
	\begin{equation}
	\label{eq:grad_flow}
	u_t=-p \,\,\,\, u|_{t=0}=f, \,\,\, p \in \partial J(u),
	\end{equation}
	where $u_t$ is the first time derivative of $u(t;x)$.
	As shown in \cite{burger2016spectral}, when the flow is initialized with an eigenfunction (that is, $\lambda f \in  \partial J(f)$), the following solution is obtained,
	\begin{equation}
	\label{eq:grad_flow_sol}
	u(t;x) = (1-\lambda t)^+ f(x),
	\end{equation}
	where $(q)^+=q$ for $q>0$ and 0 otherwise. This means that the shape $f(x)$ is spatially preserved and changes only by contrast reduction throughout time.
	An analytic solution (see \cite{benning2013ground,burger2016spectral}) can be shown for the proximal problem as well, that is, a minimization with the square $2$ norm,
	\begin{equation}
	\label{eq:ef_min_for}
	\underset{u} {\min} \,\, J(u) + \frac{\alpha}{2}\|f-u\|^2_{2}. 
	\end{equation}	
	In this case, when $f$ is an eigenfunction and $\alpha \in \mathbb{R}^+$ ($\mathbb{R}^+ = \{x \in \mathbb{R} \,| \, x \ge 0\}$) is fixed, the problem has the following solution,
	\begin{equation}
	\label{eq:ef_min_for_sol}
	u(x) = \left(1-\frac{\lambda}{\alpha}\right)^+f(x).
	\end{equation}
	In this case also, $u(x)$ preserves the spatial shape of $f(x)$ (as long as $\alpha > \lambda$).
	This was already observed by Meyer in \cite{Meyer[1]} for the case of a disk with $J$ the TV functional.
	Earlier research on \hbindex{nonlinear eigenfunctions} induced by TV, which are set indicator functions, has been referred as \emph{calibrable sets}. First aspects of this line of research can be found in the work of Bellettini et al. \cite{bellettini2002total}.
	They introduced a family of convex bounded sets $C$ with finite perimeter in $\mathbb{R}^2$ that preserve their boundary throughout the TV flow (gradient flow \eqref{eq:grad_flow} where $J$ is TV). It is shown that the indicator function of a set $C$, ${\bf 1}_C$, with perimeter $P(C)$ which admits
	\begin{equation}
	\underset{p \in \partial C}{\text{ess sup }} \kappa(p) \le \frac{P(C)}{|C|}
	\end{equation}
	is an eigenfunction, in the sense of \eqref{eq:ef_1hom}, where $u=\lambda_C {\bf 1}_C$ and
	\begin{equation}
	\label{eq:lam_tv}
	\lambda_C = \frac{P(C)}{|C|}.
	\end{equation}
	
A further generalization of  \eqref{eq:ef}, referred to as the \emph{double-nonlinear} eigenvalue problem, is formulated by introducing another bounded nonlinear operator $Q$, to have,
\begin{equation}
\label{eq:ef_double}
    T(u) = \lambda Q(u).
\end{equation}
Here $Q(u)$ may be high order polynomials or trigonometric functions. 
In physics, a variant of \eqref{eq:ef_double} is quite common, where  $T$ is a linear operator (mostly the Laplacian). For example, the one-dimensional Shcroedinger equation, 
$$-u_{xx}=\lambda(u^3-u).$$

We will address here ways also of how to solve such problems.
In the variational context, $T$ and $Q$ are two subgradient elements of different convex functionals, $J$ and $H$, thus \eqref{eq:ef_double} is rewritten as
\begin{equation}
\label{eq:ef_double_var}
    p = \lambda q,   \,\,\,\, p \in \partial J(u), \,\,\,\, q \in \partial H(u).
\end{equation}
This type of problem appears in the relaxation of the \hbindex{Cheeger} cut problem, where $J$ is TV and $H$ is $\ell^1$, see \cite{hein2010inverse,BressonSzlam2010Cheeger,feld2019rayleigh}. There are several additional algorithms which attempt to compute nonlinear eigenfunctions in some specific settings. In \cite{bozorgnia2016convergence} and \cite{bozorgnia2019approximation} algorithms for computing the smallest eigenvalue and eigenfunction of the $p$-Laplacian are proposed, along with convergence proofs. 
As part of analyzing variational networks  \cite{effland2020variational}  analyze the learned regularizers by computing their eigenfunctions. This is performed by minimizing a generalized Rayleigh quotient using accelerated gradient descent. In the process of nonlinear spectral decomposition based on gradient descent (\cite{Gilboa_spectv_SIAM_2014,burger2016spectral}), near extinction time only a single eigenfunction "survives". This idea is formalized in \cite{bungert2019computing} where eigenfunctions are computed by taking the limit at extinction time of a gradient flow. \cite{gautier2019perron,gautier2020computing} have used power-iterations to solve several nonlinear eigenpair problems. Existence and uniqueness results were obtained based on Perron-Frobenius theory. 

We will now present in detail five algorithms, coauthored by the author and colleagues, to solve various types of nonlinear eigenvalue problems. Some of the iterative algorithms can be understood as a discretization in time of a continuous nonlinear flow.

\section{Nossek-Gilboa (NG)}
This simple algorithm, presented first in \cite{nossek2018flows}, was the first of a series of algorithms, which stem from \hbindex{nonlinear flows}. These flows reach a steady-state only at eigenfunctions. Different initial conditions yield different steady-states. The goal for the (NG) algorithm is to provide a solution to the nonlinear eigenvalue problem \eqref{eq:ef_1hom}, where $J$ is an absolutely one-homogeneous functional, admitting \eqref{eq:J_1hom}. We assume a constant unit vector is in its \hbindex{null-space} (Property \ref{prop:p}). The proposed nonlinear flow is,
\begin{equation}\label{eq:flow_ng}
u_t =\frac{u}{\normd{u}}-\frac{p}{\normd{p}},\hspace{1cm} p\in\partial J(u),
\end{equation}
where $u(0)=u_0 \in X$ is an initial condition, with $\langle u_0,\mathbf{1} \rangle = 0$. 
%
The associated iterative algorithm for solving \eqref{eq:ef_1hom} is detailed in Algorithm \ref{alg:ng}.
\begin{algorithm}
\label{alg:ng}
\SetAlgoLined
\KwData{ $u_0$ with $\langle u_0,\mathbf{1} \rangle = 0$, $\Delta t \in (0,\|u_0\|_2)$, $\epsilon$.}
\KwResult{Eigenfunction and eigenvalue, $\{u^k,\lambda^k\}$, where $\lambda^k=J(u^k)/\normd{u^k}^2$.}
\textbf{Initialization:} $k\gets 0$, $u^k \gets u_0$.\\
\Repeat{$\normd{u^{k+1}-u^k}<\varepsilon$}
{
	\begin{equation}
    \label{eq:alg_ng}
	u^{k+1} = u^k + \Delta t \left( \frac{u^{k+1}}{\normd{u^k}} - \frac{p^{k+1}}{\normd{p^k}} \right),  
	\end{equation}
}
\caption{{\bf(NG).} Compute a nonlinear eigenfunction $\lambda u \in \partial J(u)$, associated with an absolutely one-homogeneous functional $J$.}
\end{algorithm}

Eq. \eqref{eq:alg_ng} is computed by solving the following convex optimization problem,
	\begin{equation}
	\label{eq:opt_ng}
	u^{k+1} = \underset{v}{\arg \min} \,\, \left\{J(v) + \frac{\normd{p^k}}{2\Delta t} \left(1 - \frac{\Delta t}{\normd{u^k}} \right) \normd{\frac{u^k}{1 - \frac{\Delta t}{\normd{u^k}}} - v}^2 \right\}.
	\end{equation}

\subsection{NG flow properties}
There are several desired properties of this flow. Although it does not emerge as a \hbindex{gradient flow} of a certain energy functional, the solution becomes smoother with time (in terms of the regularizing functional $J$). On the other hand, the $\ell^2$ norm of the solution is increasing. The main properties are summarized in the following theorem. In this case, the proof is presented, and is relatively simple to follow (it is based on \cite{nossek2018flows,aujol2018theoretical}). This allows us to get the intuition of how such flows behave. In subsequent parts, proofs are omitted and we refer the reader to the relevant papers for details, to avoid a lengthy presentation.

\begin{theorem}
Assume that there exists a solution $u$ in $W^{1,2}((0,T);X)$, $T>0$, of the flow \eqref{eq:flow_ng} . Then the following properties hold:
\begin{equation}
\frac{d}{dt} \frac12\normd{u(t)}^2 \geq 0,
 \end{equation}
moreover, we have $\langle u(t),{\bf 1}\rangle = 0 $,
and in addition,
 \begin{equation}
 \label{eq:ng_J_dec}
 \frac{d}{dt} J(u(t))\leq 0 \mbox{ for almost every $t$}.
  \end{equation}
We conclude that, $t \mapsto J(u(t))$ is non increasing for all $t \geq 0$.
\end{theorem}

\emph{Proof:}
Recalling that $\langle p,u\rangle\leq \normd{p}\normd{u}$, this flow ensures that:
\begin{equation*}
\frac{d}{dt} \frac12\normd{u(t)}^2
=\langle u,u_t\rangle=\left\langle u,\frac{u}{\normd{u}}-\frac{p}{\normd{p}}\right\rangle=\normd{u}-\frac{\langle u,p\rangle}{\normd{p}}\geq 0
\\
\end{equation*}

We can also remark that
$$\frac{d}{dt} \frac12\normd{u(t)}^2\leq \normd{u(t)}$$ so that
 $$\normd{u(t)}\leq \normd{u_0}+2t.$$

Additionally, if $u_0$ is of zero mean, Property \ref{prop:p} ensures that $u(t)$ is of zero mean, for all $t>0$.
To show \eqref{eq:ng_J_dec} we make use of Lemma 3.3 page 73 in \cite{Brezis2} (see also Lemma 4.1 in \cite{Apidopoulos}). It allows us to use the "chain rule for differentiation". Let us first recall this lemma.
\begin{lemma}[Brezis '73]
\label{lem:Brezis}
Let $T>0$ and $F$ be a convex, lower semi-continuous, proper function and $v\in W^{1,2}((0,T);X).$ Let also $h\in L^2((0,T);X)$, such that $h\in \partial F(v(t))$ a.e. in $(0,T)$. Then the function $F \circ v:[0,T]\to \mathbb{R}$ is absolutely continuous in $[0,T]$ with 
$$ \frac{d}{dt}\left(F(v(t))\right) = \langle z,v_t \rangle,  \,\,\, \forall z \in \partial F(v(t)) \,\,\, a.e. \textrm{ in } (0,T).$$ 
\end{lemma}
From Lemma \ref{lem:Brezis}, if $u$ is in $W^{1,2}((0,T);X)$, we get that $J(u(t))$ is absolutely continuous in $[0,T]$ with
 \begin{equation*}
 \frac{d}{dt} J(u(t))
 =\langle p,u_t\rangle=\left\langle p,\frac{u}{\normd{u}}-\frac{p}{\normd{p}}\right\rangle=\frac{\langle u,p\rangle}{\normd{u}}-\normd{p}\leq 0.
 \end{equation*}
This inequality holds for almost every $t$, and since  $t \mapsto J(u(t))$ is an absolutely continuous function, we deduce that it is a non increasing function. 
$\square$

The flow \eqref{eq:flow_ng} converges iff  $u_t=0$ so that
$$p=\frac{\normd{p}}{\normd{u}}u\in\partial J(u)\Rightarrow p=\frac{J(u)}{\normd{u}^2}u$$
and $u$ is an eigenfunction of $J$ with eigenvalue $\lambda=\frac{J(u)}{\normd{u}^2}$.

\subsection{NG iterations algorithm properties}
The iterations in Algorithm \ref{alg:ng} can be viewed as a \hbindex{semi-implicit scheme} of the flow \eqref{eq:flow_ng}.
The properties of the discrete flow are similar in nature to those of the continuous flow (but not precisely the same). They are summarized in the following theorem (details are given in \cite{nossek2018flows}).
\begin{theorem}
\label{th:ng_flow_discrete}
The solution $u^k$ of the discrete flow \eqref{eq:alg_ng} of Algorithm \ref{alg:ng} has the following properties:
\begin{itemize}[labelindent=1.0em,labelsep=0.5cm,leftmargin=*]
			\item[(i)]  $\langle  u^k,\mathbf{1} \rangle =0. $
			\item[(ii)] $\|p^{k+1}\|_2\le \| p^k\|_2. $ 
			\item[(iii)] $ \| u^{k+1}\|_2 \ge \|u^k \|_2. $ 
            \item[(iv)] $\frac{J(u^{k+1})}{\normd{u^{k+1}}}  \le \frac{J(u^k)}{\normd{u^{k}}} . $ 
			\item[(v)] A sufficient and necessary condition for steady-state $u^{k+1} = u^k$ holds if $u^k$ is an eigenfunction, admitting \eqref{eq:ef_1hom}.
\end{itemize}
\end{theorem}

\section{Aujol-Gilboa-Papadakis (AGP)}
In \cite{aujol2018theoretical} the authors proposed a generalized flow for solving \eqref{eq:ef_1hom}, which is more stable than (NG) and can be better analyzed theoretically. The general flow, for $\alpha\in[0;1]$, is,
\begin{equation}
\label{eq:flow_agp_alpha}
u_t=\left(\frac{J(u)}{\normd{u}^2}\right)^\alpha u-\left(\frac{J(u)}{\normd{p}^2}\right)^{1-\alpha}p,\hspace{1cm} p\in\partial J(u),
\end{equation}
with $u(0)=u_0 \in X$, $\langle u_0,\mathbf{1} \rangle = 0$.
Notice that for $\alpha=1/2$, we retrieve the (NG) flow,  \eqref{eq:flow_ng}, up to a normalization with $J^{1/2}(u)$.
For the case $\alpha=1$ the flow becomes  
\begin{equation}
\label{eq:flow_agp}
u_t=\left(\frac{J(u)}{\normd{u}^2}\right) u - p,\hspace{1cm} p\in\partial J(u).
\end{equation}
In this case, there is no term with $\normd{p}$ in the denominator and the analysis simplifies. \hbindex{Uniqueness} of the flow and \hbindex{convergence} of the iterative algorithm are established.

For the case $\alpha=1$ we get that the  $\ell^2$ norm is fixed in time. This allows us to have a unit norm throughout the evolution. In the discrete iterations, however, an additional normalization step is required to maintain this property.
Given any input $f$, to obtain a valid initial condition $u_0$, we first subtract the mean and then normalize by the $\ell^2$ norm .
The associated iterative algorithm, $\alpha=1$, for solving \eqref{eq:ef_1hom} is detailed in Algorithm \ref{alg:agp}.
\begin{algorithm}
\label{alg:agp}
\SetAlgoLined
\KwData{ $u_0$ with $\langle u_0,\mathbf{1} \rangle = 0$, $\normd{u_0}=1$, $\Delta t \in (0,\normd{u_0}^2/J(u_0))$, $\epsilon$.}
\KwResult{Eigenfunction and eigenvalue, $\{u^k,\lambda^k\}$, where $\lambda^k=J(u^k)/\normd{u^k}^2$.}
\textbf{Initialization:} $k\gets 0$, $u^k \gets u_0$.\\
\Repeat{$\normd{u^{k+1}-u^k}<\varepsilon$}
{
	\begin{equation}
    \label{eq:alg_agp}
    \begin{split}
	u^{k+1/2} &= u^k + \Delta t \left( \frac{J(u^k)u^{k+1/2}}{\normd{u^k}^2} - p^{k+1/2} \right), \\
	u^{k+1} &= \frac{u^{k+1/2}}{\normd{u^{k+1/2}}}.
	\end{split}
	\end{equation}
}
\caption{{\bf(AGP).} Compute a nonlinear eigenfunction $\lambda u \in \partial J(u)$, associated with an absolutely one-homogeneous functional $J$.}
\end{algorithm}


The term $u^{k+1/2}$ in Eq. \eqref{eq:alg_agp} is computed by solving,
	\begin{equation}
	\label{eq:opt_agp}
	u^{k+1/2} = \underset{v}{\arg \min} \,\, \left\{J(v) +
	\frac{1}{2 \Delta t}\normd{v-u^k}^2 
	-\frac{J(u^k)}{2\normd{u^k}^2}\normd{v}^2 
    \right\}.
	\end{equation}
There is a unique minimizer $v$ for any time step $\Delta t$ which is in the range specified above.

\subsection{AGP flow properties}

\begin{theorem}\label{prop_alpha}
For $u_0$ of zero mean and  $\forall\alpha\in[0;1]$, if $u$ is in $ W^{1,2}((0,T);X)$, then the trajectory $u(t)$ of the flow \eqref{eq:flow_agp_alpha} satisfies the following properties: 
\begin{itemize}[labelindent=1.0em,labelsep=0.5cm,leftmargin=*]
\item [(i)]$\langle u(t),\mathbf{1}\rangle=0.$
\item [(ii)] $\frac{d}{dt}  J(u(t))\leq 0$ for almost every $t$. Moreover, $t \mapsto J(u(t))$ is non increasing. If  $\alpha=0$, we have for almost every $t$ that   $\frac{d}{dt}  J(u(t))= 0$ and $t \mapsto J(u(t))$ is constant.
\item [(iii)]$\frac{d}{dt} \normd{u(t)}\geq 0$ and  $\frac{d}{dt}  \normd{u(t)}= 0$ for $\alpha=1$.
\item[(iv)] If the flow converges to $u^*$,  we have $p^*=J^{2\alpha-1}(u^*)\frac{\normd{p^*}^{2(1-\alpha)}}{\normd{u^*}^{2\alpha}}u^*\in\partial J(u^*)$ so that $u^*$ is an eigenfunction.\\
\end{itemize} 
\end{theorem}

\noindent {\bf Uniqueness.} For the case $\alpha=1$, one can establish
uniqueness of the flow \eqref{eq:flow_agp}, under mild conditions. 
\begin{theorem}
Let  $u$ and $v$ be two solutions of \eqref{eq:flow_agp} in $ W^{1,2}((0,T);X)$ with respective initial condition $u_0$ and $v_0$, such that $J(u_0) < + \infty$ and  $J(v_0) < + \infty$, with $\|u_0\|_2=\|v_0\|_2=1$.
Then we have:
\begin{equation} \label{res_prop}
\frac{d}{dt} \left( \frac{1}{2} \|u-v\|_2^2 \right) \leq
\frac{J(u)+J(v)}{2}  \|u-v\|_2^2 .
\end{equation}
\end{theorem}
By the fact that $J(u)$ is decreasing and using \hbindex{Gronwall lemma} we obtain
\begin{equation} \label{uniqueness}
\|u-v\|_2^2
\leq
\|u_0-v_0\|_2^2
\exp \left( (J(u_0)+J(v_0)(t-t_0) \right).
\end{equation}

\subsection{AGP iterations algorithm properties}
The iterations in Algorithm \ref{alg:agp} can be viewed as a semi-implicit scheme of the flow \eqref{eq:flow_agp}. The algorithm's properties are detailed below. 
\begin{theorem}
\label{th:agp_flow_discrete}
Let $u_0$ in $X$, and the sequence $u_k$ defined by \eqref{eq:alg_agp}. Then the sequences $J(u_k)$ and $\|p_k\|_2$ are non increasing, $\|u_k\|_2 =\|u_0\|_2$ for all $k$, and $u_{k+1}-u_k \to 0$.
\end{theorem}

\noindent {\bf Convergence.} Finally, it is shown that Algorithm \ref{alg:agp} converges to an eigenfunction.
\begin{theorem} \label{th_conv_finite}
Let $u_0$ be in $X$, and the sequence $u_k$ be defined by \eqref{eq:alg_agp}. 
There exist some $u$ and $p$ in $X$ such that, up to a subsequence, $u_k$ converges to $u$ in $X$ and $p_k$ converges to $p$ in $X$, with $p \in \partial J(u)$, and $J(u_k)$ converges to $J(u)$.
Moreover, $u$ is a nonlinear eigenfunction, in the sense of \eqref{eq:ef_1hom}.
\end{theorem}

\section{Feld-Aujol-Gilboa-Papadakis (FAGP)}
In \cite{feld2019rayleigh} the aim is to solve the problem \eqref{eq:ef_double_var} for the case when $J$ and $H$ are both absolutely one-homogeneous functionals.
Let us consider the generalized nonlinear  \hbindex{Rayleigh quotient} 
\begin{equation}\label{eq:rayleigh}
R(u):=\frac{J(u)}{H(u)}.
\end{equation}
In an analogue to the linear case, eigenfunctions in the sense of \eqref{eq:ef_double_var} are  critical points of \eqref{eq:rayleigh}. In segmentation, classification and clustering, often we seek eigenfunctions with the least (strictly positive) eigenvalue. Thus, excluding the null-space of $J$ and $H$, we seek to minimize the Rayleigh quotient \eqref{eq:rayleigh}.
%
A classical way to reach a local minimizer of $R(u)$ is by using a gradient descent flow, 
$$ u_t = -\nabla R(u).$$
Taking the variational derivative of $R(u)$, with $q\in\partial H(u), p\in\partial J(u)$, the gradient descent flow is, 
\begin{equation}\label{gd_flow}
u_t= \frac{J(u)q - H(u)p}{H^2(u)}.  
\end{equation}
The flow can also be written as, 
$$ u_t = \frac{R(u)q-p}{H(u)}. $$
This flow is hard to analyze theoretically, mainly due to the division by $H(u)$. 
Therefore, \cite{feld2019rayleigh} proposed the following flow to minimize $R(u)$,
\begin{equation}\label{main_flow}
u_t= R(u) q - p.
\end{equation}
This is essentially a gradient-descent type flow, without the division by $H(u)$, which can be interpreted as a dynamic rescaling of the time parameter. The flow reduces monotonically the quotient $R(u)$ and the steady state admits the nonlinear eigenvalue problem \eqref{eq:ef_double_var}. 

A second flow is proposed, that minimizes the log of the Rayleigh quotient,
$$ u_t = -\nabla (\log R(u)), $$
which can be written as, 
\begin{equation}\label{log_flow}
u_t= \frac{q}{H(u)}- \frac{p}{J(u)}.
\end{equation}
This is motivated by a widely used practice of using the log of a function involving multiplicative expressions. It is commonly employed in statistics and machine learning algorithms, such as maximum likelihood estimation and policy learning.
The flow is essentially a time rescaling of \eqref{main_flow} by $1/J(u)$. We note that it is not in the form of Brezis Lemma \ref{lem:Brezis} and therefore is harder to analyze. We will not focus on this flow here. It is worth mentioning, however, that  
in the context of the Cheeger cut problem, we found out that numerically it is very stable and highly resilient to the choice of the discrete time step. Thus a large time step can be chosen, which speeds up numerical convergence (see details in \cite{feld2019rayleigh}).

The algorithm is based on the following semi-explicit scheme of the flow,
\begin{equation}\label{flow_main1}
\left\{\begin{array}{ll}
(u^{k+1/2}-u^k)/\Delta t=R(u^k)q_k-p_{k+1/2},\ \ \ q_k\in\partial H(u^k),\,\,p_{k+1/2}\in\partial J(u^{k+1/2})\\
u^{k+1}=u^{k+1/2}/\normd{u^{k+1/2}}.
\end{array}
\right.
\end{equation}
This scheme is associated with the minimization of a convex functional,
\begin{equation}\label{general_func}
u^{k+1/2}=\uargmin{u \in X}F(u):=\frac1{2 \Delta  t}\normd{u-u^k}^2- R(u^k)\left\langle q_k,u\right\rangle+J(u),
\end{equation}
where $u^{k+1/2}$ being a minimizer of $F$ implies that there exist $p_{k+1/2}\in\partial J(u^{k+1/2})$ such that  
$$\frac{1}{d t}(u^{k+1/2}-u^k)- R(u^k) q_k+p_{k+1/2}=0.$$
This leads directly to Algorithm \ref{alg:fagp}.

\begin{algorithm}\label{alg:fagp}
\SetAlgoLined
\KwData{ $u_0$ with $\langle u_0,\mathbf{1} \rangle = 0$, $\normd{u_0}=1$, $\Delta t>0$, $\epsilon>0$.}
\KwResult{Local Minimizer $u$ of the Rayleigh quotient $R=J/H$.}
\textbf{Initialization:} $k\gets 0$, $u^k \gets u_0$.\\
\Repeat{$\normd{u^{k+1}-u^k}<\varepsilon$}
{$u^{k+1/2}=\uargmin{u \in X}F(u):=\frac1{2 \Delta t}\normd{u-u^k}^2- R(u^k)\left\langle q^k,u\right\rangle+J(u).$\\
$u^{k+1}=u^{k+1/2}/\normd{u^{k+1/2}}$\\
}
\textbf{end while}

\caption{(FAGP).Rayleigh quotient minimization of absolutely one-homogeneous functionals}

\end{algorithm}
\vspace{2mm}

\begin{rem}\label{normalization_step}
Notice that since J and H are absolutely one-homogeneous their subgradients do not change by the normalization step of the flow, i.e $q_{k+1}=q_{k+1/2}$ and $p_{k+1}=p_{k+1/2}$. We also have $R(u^{k+1})=R(u^{k+1/2})$ as a quotient of two one-homogeneous functionals.
\end{rem}
 The sequence $u^k$ of Algorithm \ref{alg:fagp} satisfies the following properties:
\begin{enumerate}
\item  $1=\normd{ u^{k}}^2\leq\langle  u^{k+1/2},u^k\rangle \leq \normd{u^{k+1/2}}^2$.  
\item \label{point2} $\normd{u^{k+1}-u^k}\leq \normd{u^{k+1/2}-u^k}$.
\item Monotonicity: $R(u^{k+1})\leq R(u^k)$.
\item Compactness: $\normd{u^{k+1}-u^k}^2 \to 0$.
\end{enumerate}

\noindent {\bf Convergence.} It is shown that Algorithm \ref{alg:fagp} converges to a (double nonlinear) eigenfunction, in the sense of \eqref{eq:ef_double_var}.
\begin{theorem}
[Convergence]\label{discrete_convergence}
Let $u_0$ in $X$ and $u^k$ is computed by Algorithm \ref{alg:fagp}. Then there exist $u$, $p$ and $q$ in $X$ such that up to a subsequence $u^k \to u$, $p_{k+1/2} \to p$, $q_k \to q$,  $\normd{u}=1$, and
\begin{equation}
p=R(u)q,\,\,\,\,q\in\partial H(u),\,\,p \in\partial J( u).
\end{equation}
\end{theorem}

Further relations to calibrable sets and variants of Algorithm \ref{alg:fagp} for Cheeger cut minimization on graphs are provided in detail in \cite{feld2019rayleigh}.

\section{Cohen-Gilboa (CG)}
Nonlinear eigenvalue problems emerge naturally also in physical modeling of nonlinear phenomena in fields such as photo-electronics and quantum physics. 
In 1895	Korteweg-de Vries formulated a mathematical model of waves on shallow water surfaces which were previously described by Russell. 
The \hbindex{KdV equation}, as expressed in \cite{zabusky1965interaction}, is,
\begin{equation*}
	u_t + uu_x+\delta^2u_{xxx}=0,
\end{equation*}
with $\delta$ a small real scalar. Reformulating this expression for a stationary wave yields,
\begin{equation}\label{eq:KdV}
	-u_{XX}=\lambda\left(-cu+\frac{u^2}{2}\right),
\end{equation}
where $c$ is the wave velocity, $X=x-ct$ and $\lambda=\delta^{-2}$. Naturally, $\lambda$ can be understood as an eigenvalue.
The solution to this equation models well a family of solitary waves  referred to as \hbindex{solitons}.
In this specific case one can obtain an analytic solution,
\begin{equation*}
	u(X) = 3c\cdot\textrm{sech} ^2\left(
    \frac{\sqrt[]{c\cdot\lambda}X}{2}
    \right).
\end{equation*}
In recent decades there has been a growing research concerning nonlinear physical models, where more complex nonlinear eigenvalue problems emerge, such as the two-dimensional nonlinear Schroedinger equation, 
\begin{equation}
\label{eq:nls}
	u_{xx}+u_{yy}-V_0\left(sin^2x +sin^2y\right)u+\sigma |u|^2u=-\mu u.
\end{equation}
In \cite{cohen2018energy} a method for solving such problems was proposed, following the flows of \cite{nossek2018flows} and \cite{aujol2018theoretical}. The basic formulation was to solve the (double) nonlinear problem,
\begin{equation}
\label{eq:cg_ef}
    T(u) = \lambda Q(u),
\end{equation}
where $T(u)\in \partial J(u)$, $J(u)$ is a convex, proper, lsc regularizing functional and $Q(u)$ is a bounded nonlinear operator, with both $T,Q \in L^2(\Omega)$ . The following flow is a natural generalization of \cite{nossek2018flows}, 
\begin{equation}\label{eq:mainFlow}
	u_t(t) = M(u(t)), \,\,\,\,\,\, u(t=0)=u_0,
\end{equation}
where 
\begin{equation}
\label{eq:Mu}
    M(u) = s\frac{Q(u)}{\normd{Q(u)}}-\frac{T(u)}{\normd{T(u)}},
\end{equation}
and $s = \textrm{sign}(\langle {Q(u)},{T(u)}\rangle)$.
It can be shown that $\frac{d}{dt}J(t) \le 0$ a.e. for $t\in (0,\infty)$ and that a steady state admits the nonlinear eigenvalue problem \eqref{eq:cg_ef}. 

A problem arises here, where one can reach the null-space of $J$, thus yielding degenerate solutions with eigenvalues $\lambda=0$. This did not happen in previous algorithms, which ensured $u$ to be of zero mean and unit norm (or increasing norm with time in \cite{nossek2018flows}). This prevented the case where $u$ can be a constant function. For \eqref{eq:cg_ef}, however, these assumptions do not necessarily hold, moreover we do not control $u$ directly. Such flows tend to find smoother solutions with low eigenvalues, thus reaching a very smooth degenerate solution is not only a theoretical problem, but a phenomenon which is actually encountered in numerical experiments. Thus, one needs to "push" the evolution "away" from degenerate solutions. This is formulated in general by defining a subspace which does not include all eigenfunctions with zero eigenvalues. We would like our flow to always stay in that subspace. An additional term is added to the flow, which directs  it toward this subspace. Let us explain it in more details for the case where $J$ is the Dirichlet energy, hence $T(u)=-\Delta u$. We thus want to solve,
\begin{equation}
\label{eq:cg_ef_lap}
    -\Delta u = \lambda Q(u).
\end{equation}
This is an eigenvalue problem with left-sided linear operator and right-sided nonlinear operator (common in physics).
For Neumann boundary conditions the null space of $J$ is the space of constant functions. Therefore, the following energy is defined,
\begin{equation}
E(u) = \frac{1}{2}\langle{Q(u)},{1}\rangle^2,
\end{equation}
with 
$$ \partial E = \langle{Q(u)},{1}\rangle \partial Q, $$
and $\partial Q$ is the variational derivative of $\langle{Q(u)},{1}\rangle$.
We would like $E(u)=0$ at steady-state to ensure we obtain a meaningful solution. A variant of a gradient descent with respect to $E$ is defined by,
\begin{equation}\label{eq:compFlow}
	u_t=C(u)
\end{equation}
where
\begin{equation}
\label{eq:Cu}
	C(u) = -\partial_{u} E + \frac{\langle {\partial_{u} E},{T(u)}\rangle}{\normd{T(u)}^2}T(u).
\end{equation}
It ensures one decreases $E$ while not increasing $J$. We call this the complementary flow.
Let us compute the time derivatives of $J$ and $E$:
\begin{equation}\label{eq:compFlowProofJ}
\begin{split}
\frac{d}{dt}J(u)=&\langle{T(u)},{u_t}\rangle=\langle{T(u)},{C(u)}\rangle\\
=&\langle{T(u)},{-\partial_{u} E + \frac{\langle{\partial_{u} E},{T(u)}\rangle}{\normd{T(u)}^2}T(u)}\rangle=0.
\end{split}
\end{equation}
For $E$ we have,
\begin{equation}\label{eq:compFlowProofE}
\begin{split}
\frac{d}{dt}E(u)=&\langle{\partial_uE},{u_t}\rangle=\langle{\partial_{u} E},{C(u)}\rangle\\
=& -\normd{\partial_{u} E}^2 + \frac{\langle{\partial_{u} E},{T(u)}\rangle^2}{\normd{T(u)}^2}\leq 0,
\end{split}
\end{equation}
where the last inequality follows Cauchy-Schwarz. We thus can merge the main flow \eqref{eq:mainFlow} and the complementary one \eqref{eq:compFlow}, with some weight parameter $\alpha$ to obtain the final flow,
\begin{equation}\label{eq:combinedFlow}
	u_t = M(u)+\alpha C(u),
\end{equation}
where $\alpha\in\mathbb{R}_+$ and $M(u)$ and $C(u)$ are defined in \eqref{eq:Mu}, and \eqref{eq:Cu}, respectively. This combined flow admits $(d/dt)J(u) \le 0 $ and $(d/dt)E(u)\le 0$ (for $\alpha$ large enough).
Numerically, iterations which follow this flow are provided in 
\cite{cohen2018energy}, using the following adaptive time step for the main flow,
\begin{equation}\label{eq:dtM}
dt_M = 2\frac{\langle{\Delta u^{k}},{M(u^k)}\rangle}{\normd{\nabla M(u^k)}^2},
\end{equation}
and an adaptive step size for the complementary flow, 
\begin{equation}\label{eq:dtC}
dt_C = -\frac{E(u^{k+\frac{1}{2}})}{\langle{\partial E(u^{k+\frac{1}{2}})},{C(u^{k+\frac{1}{2}})}}\rangle.
\end{equation}
The choice of $dt_C$ was such that it approximates in a single step $E(u)\approx 0$, within a first Taylor approximation.
The numerical algorithm, a dissipating flow with respect to the energy term $J$ (ensured to be non-increasing), is shown in Algorithm \ref{alg:cg}. Since it is basically an explicit scheme with carefully chosen time-steps, each iteration requires a low computational effort.

\begin{algorithm}
\label{alg:cg}
\SetAlgoLined
\KwData{ $u_0$, $Q(u)$, $\epsilon>0$.}
\KwResult{Eigenfunction and eigenvalue, $\{u^k,\lambda^k\}$, where $\lambda^k= \langle T(u),u \rangle / \langle Q(u),u \rangle$.}
\textbf{Initialization:} $k\gets 1$, $u^k \gets u_0$, $T(u)=-\Delta u$.\\
Set $dt_C(u_0)$ according to \eqref{eq:dtC}.\\
$u^1 \gets u^0+ dt_C(u_0) \cdot C(u^0)$.\\
\Repeat{$\normd{u^{k+1}-u^k}<\varepsilon$}
{
    Set $dt_M$ according to \eqref{eq:dtM} and $M(u^k)$ according to \eqref{eq:Mu}.\\
        	$u^{k+\frac{1}{2}} \gets u^{k}+dt_M\cdot M(u^k)$.\\
			Set $dt_C$ according to \eqref{eq:dtC} and $C(u^{k+\frac{1}{2}})$ according to \eqref{eq:Cu}.\\
        	$u^{k+1} \gets u^{k+\frac{1}{2}}+dt_C  \cdot C(u^{k+\frac{1}{2}}). $
}
\caption{{\bf(CG).} Nonlinear eigenpair generation for the Laplacian problem: $-\Delta u = \lambda Q(u)$.}
\end{algorithm}

\section{Bungert-Hait-Papadakis-Gilboa (BHPG)}
The last algorithm presented here is related to very general and complex nonlinear operators, which often cannot be expressed analytically. In \cite{hait2019numeric} and  \cite{bungert2020powermethod} the operators considered were nonlinear \hbindex{denoisers}, which can be based on classical algorithms or on \hbindex{deep neural networks}. 

The setting is as follows. Let $T:\H\to\H$ be a generic (nonlinear) operator on a real Hilbert space $\H$ with norm $\norm{\cdot}$.
In the case of a neural network one typically has $\H=\R^n$, equipped with the Euclidean norm.
We aim at solving the nonlinear eigenproblem  \eqref{eq:ef},  
$$    T(u)=\lambda u, $$
where $u\in\H$ and $\lambda\in\R$ denote the eigenvector and eigenvalue, respectively.
Since the operator assumed here is very general and is not based on any energy functional, one needs to resort to a very simple iterative process, which does not involve any minimization. Such a simple algorithm exists for the linear case, the \hbindex{power method}.

Linear power method is a simple classical algorithm for solving linear eigenvalue problems $Lu = \lambda u$, where $u\in \R^n$ is a vector and $L\in \mathbb{R}^{n \times n}$ is a diagonalizable matrix. Given some initial condition $u_0$, $k\gets 0$, $u^k \gets u_0$, the following process is iterated until convergence,
\begin{equation}
    \label{eq:pm_lin}
    u^{k+1} \gets \frac{Lu^k}{\|Lu^k\|_2}, \quad \,\,k \gets k+1.
\end{equation}
Under mild conditions, it is known to converge to the eigenvector with the largest eigenvalue, although convergence is slow.
A straightforward analog of this process for the nonlinear case, having an operator $T(u)$, is to initialize similarly and to iterated until convergence,
\begin{equation}
    \label{eq:pm_nonlin}
    u^{k+1} \gets \frac{T(u^k)}{\|T(u^k)\|_2}, \quad \,\,k \gets k+1.
\end{equation}
One can analyze this process more easily in a restricted nonlinear case, where $J$ is an absolutely one-homogeneous functional, based on a proximal operator of $J$,
\begin{align}\label{eq:prox}
    \prox_\alpha^J(u):=\argmin_{v\in\H}\frac{1}{2}\norm{v-u}^2+\alpha J(v),
\end{align}
where $u\in\H$ and $\alpha > 0$ denotes the regularization parameter. The operator is a classical variational denoiser, 
\begin{equation}\label{eq:Tu_prox}
    T(u) = \prox_{\alpha}^J(u),
\end{equation}
which for $J=TV$ coincides with the ROF denoising model (\cite{rof92}). 
In \cite{bungert2020powermethod} it was shown that the process is well defined for a range of parameters $\alpha$, that the energy is decreasing, $J(u^{k+1}) \le J(u^k)$, along with a full proof of convergence to a nonlinear eigenvector, in the sense of \eqref{eq:ef}.

For more complex nonlinear operators, however, certain modifications are required. A critical issue is the range of the operator. Unlike linear or homogeneous operators, general nonlinear operators often are expected to perform only in a certain range. This is certainly true in neural-networks, where the range is dictated implicitly by the range of the images in the training set. Thus normalization by the norm, as in \eqref{eq:pm_nonlin}, can drastically change the range of $u^k$ and cause unexpected behavior of the operator. Furthermore, the mean value of $u^k$ is a significant factor. For denoisers, we often expect that a denoising operation does not change the mean value of the input image, that is   
\begin{equation}
    \label{eq:meanTu}
    \langle T(u),1 \rangle = \langle u,1 \rangle.
\end{equation}
It can be shown that for any vector $u\neq 0$ with non-negative entries and a denoiser $T$ admitting \eqref{eq:meanTu}, if $u$ is an eigenvector then $\lambda=1$.
Another issue is the invariance to a constant shift in illumination. We expect the behavior of $T$ to be invariant to a small global shift in image values. That is, $T(u+c) = T(u) + c$, for any $c \in \mathbb{R}$, such that $(u+c) \in \H$.

We thus relax the basic eigenproblem \eqref{eq:ef} as follows:
\begin{align}\label{Eq:new_eigenproblem}
  {T(u)-\overline{T(u)}}= \lambda ({u-\overline{u}}),
\end{align}  
where $\lambda\in\mathbb{R}$, $\bar u=\langle 1,u \rangle /|\Omega|$ is the mean value of $u$ over the image domain $\Omega$. 
Note that now (relaxed) eigenvectors, admitting  \eqref{Eq:new_eigenproblem}, can have any eigenvalue, keeping the assumptions on $T$ stated above. In addition, if $u$ is an eigenvector, so is $u+c$, as expected for operators with invariance to global value shifts.
A suitable Rayleigh quotient, associated with the relaxed eigenvalue problem \eqref{Eq:new_eigenproblem}, is,
\begin{align}
\label{eq:Rayleigh_new}
    R^\dagger(u)=\frac{\langle u-\overline{u},T(u)-\overline{T(u)}\rangle}{\normd{u-\overline{u}}^2},
\end{align}
which still has the property that $\lambda=R^\dagger(u)$ whenever $u$ fulfills \eqref{Eq:new_eigenproblem}.
The modified nonlinear power method is detailed in Algorithm \ref{alg:bhpg}, aiming at computing a relaxed eigenvector \eqref{Eq:new_eigenproblem} by explicitly handling the mean value and keeping the norm of the initial condition.
We found this adaptation to perform well on denoising networks.

\begin{algorithm}
\label{alg:bhpg}
\caption{(BHPG): Nonlinear power method for non-homogeneous operators.}
\SetAlgoLined
\KwData{ $u_0$, $\epsilon>0$.}
\KwResult{Relaxed eigenpair $(u^*, \lambda^*)$ in the sense of \eqref{Eq:new_eigenproblem}, where $u^*=u^k$, $\lambda^*=R^\dagger(u^*)$, with $R^\dagger$ defined in \eqref{eq:Rayleigh_new}.}
\textbf{Initialization:} $k\gets 0$, $u^k \gets u_0$.\\
\Repeat{$\normd{u^{k+1}-u^k}<\varepsilon$}
{
    $u^{k+1} \gets T(u^k)$.\\
    $u^{k+1} \gets u^{k+1}-\overline{u^{k+1}}$.\\
    $u^{k+1} \gets \frac{ u^{k+1}}{\|u^{k+1}\|}\|u_0-\overline{u_0} \|$.\\
    $u^{k+1} \gets u^{k+1} + \overline{u^k}, \,\,\,k \gets k+1.$
}
\end{algorithm}

\section{Evaluation and Examples}
We present here several results of the algorithms presented earlier. First we discuss how the numerical solutions can be evaluated. Then we show several numerical examples related to  image processing, learning and physics.

\begin{figure}[htb]
	\centering
	\begin{tabular}{cc}
   \includegraphics[width=0.35\textwidth]{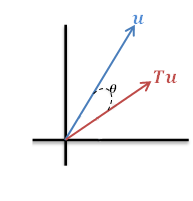}&
\includegraphics[width=0.35\textwidth]{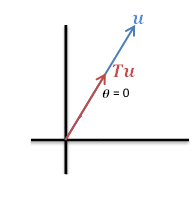}
\end{tabular}
	\caption{Global measure $\theta$, \eqref{eq:Theta}. Measures the angle between $u$ and $T(u)$. For $theta = 0$ we have a precise eigenfunction (also for $180$ degrees, negative eigenvalues).}
	\label{fig:theta}
\end{figure} 

\begin{figure}[htb]
	\centering
	\begin{tabular}{cc}
 \includegraphics[width=0.35\textwidth,height=0.2\textwidth]{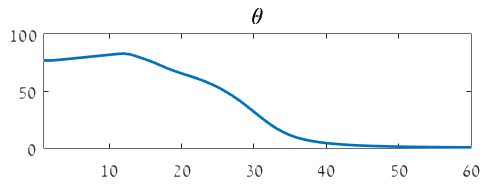}
 \includegraphics[width=0.35\textwidth,height=0.2\textwidth]{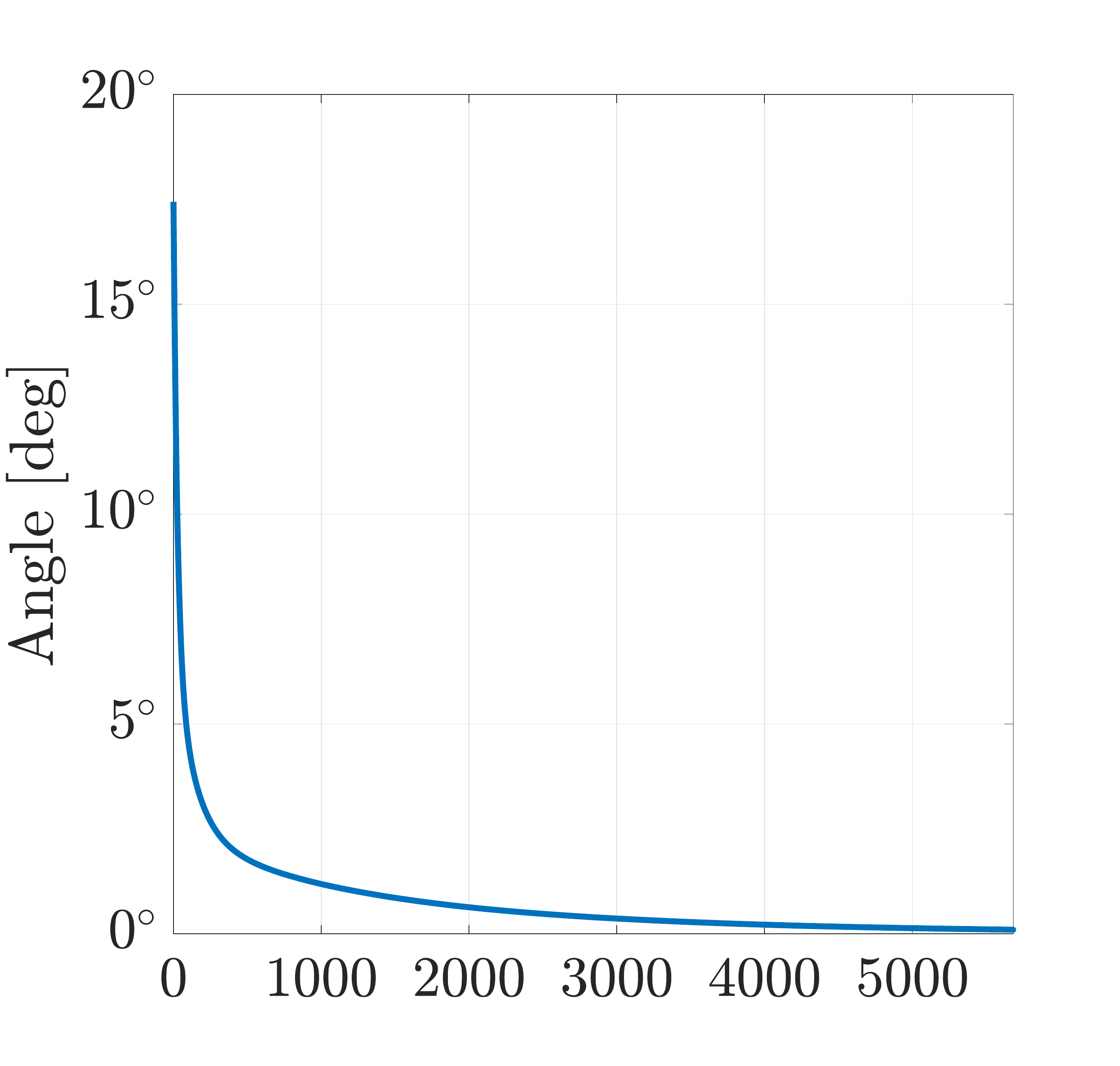}
	\end{tabular}
	\caption{$\theta$ (degrees) as a function of iterations, for (NG) flow, $J=TGV$, and for (CG) flow, Nonlinear Schrodinger equation. Taken from \cite{nossek2018flows} and \cite{cohen2018energy}.}
	\label{fig:theta_time}
\end{figure} 

\begin{figure}[htb]
	\centering
	\begin{tabular}{ccc}
\includegraphics[width=0.30\textwidth]{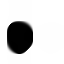}&
\includegraphics[width=0.30\textwidth]{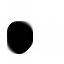}&
\includegraphics[width=0.30\textwidth]{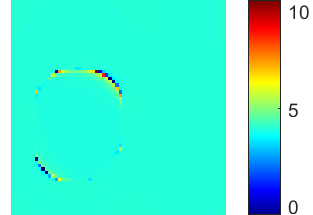}\\
TV: $u \qquad$& $T(u)=p \qquad$ & $T(u)/u$\\
\includegraphics[width=0.30\textwidth]{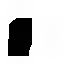}&
\includegraphics[width=0.30\textwidth]{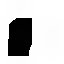}&
\includegraphics[width=0.30\textwidth]{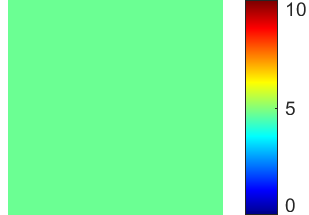} \\
ATV: $u \qquad$& $T(u)=p \qquad$ & $T(u)/u$\\
	\end{tabular}
	\caption{Local measure $\Lambda(x)=T(u)/u$. At convergence $T(u)=\lambda u$, thus for any $u\ne0$ we can examine the ratio  $\Lambda(x)$, which should be a constant function of value $\lambda$, $\forall x$. Top row, algorithm did not fully converge yet, $u$ is close to an eigenfunction for isotropic TV, the ratio (right) exposes areas where there is deviation from a constant. Bottom row, a converged eigenfunction for anisotropic TV. The ratio image is constant, up to numerical precision. Taken from \cite{aujol2018theoretical}.  }
	\label{fig:pu_ratio}
\end{figure}

\subsection{Global and local measures}
Since there is often no ground truth or analytic solutions for nonlinear eigenvalue problems, we need to find alternative ways to determine whether the algorithm converged to an eigenfunction. Often exact convergence is very slow, thus knowing that you approximately reached an eigenfunction numerically may also speed up the algorithm and serve as a good stopping criterion for the iterative process.

 One general formulation for any operator $T$, is by the angle (see~\cite{nossek2018flows}).
For eigenvectors, vectors $u$ and $T(u)$ are collinear. Thus their respective angle is either~$0$ (for positive eigenvalues) or~$\pi$ (for negative eigenvalues). Since both $u$ and $T(u)$ are real, eigenvalues are also real. Thus, the angle is a simple scalar measure that quantifies how close $u$ and $T(u)$ are to collinearity. We define the angle $\theta$ between $u$ and $T(u)$ by
\begin{equation}
\label{eq:Theta}
    \cos (\theta) = \frac{\langle u,T(u)\rangle}{\|u\|\| T(u)\|}.
\end{equation}
See Fig. \ref{fig:theta} for an illustration of $\theta$.
In most cases discussed here we have positive eigenvalues, thus we aim to reach an angle close to 0.
In Fig. \ref{fig:theta_time} we show two examples of the behavior of theta over time for (NG) and (CG) algorithms. Note that $\theta$ may not be monotonic and may increase in some time range.
The angle $\theta$ is a good global measure. In the iterative algorithms, it can be used as a stopping criteria. Instead of requiring $\normd{u^{k+1}-u^k}<\varepsilon$ one can require reaching a small enough theta $\theta < \theta_{thres}$. In our studies we often regard a function with $\theta < \pi/360$ ($\frac{1}{2}$ degree) as a numerical eigenfunction. 

One may also like to have a local measure. Usually there is no precise pointwise convergence of $(T(u))(x)=\lambda u(x)$, $\forall x$. A good way to see how spatially  the function is close to an eigenfunction is by examining the ratio
$$ \Lambda(x) = \frac{T(u)}{u},  \quad \forall u(x)\ne 0. $$
At full convergence we should have $\Lambda(x) \equiv \lambda$. The deviation map from a constant function reveals the areas where the numerical approximation is less accurate. To avoid dividing by values close to 0, one may compute this map only for $u(x) >  \delta$, where $\delta$ is a small constant.
In Fig. \ref{fig:pu_ratio} we show two examples of this ratio, when one obtains a function close (but not precisely) an eigenfunction and for a case with full convergence.


\begin{figure}[htb]
	\centering
	\begin{tabular}{cccc}
	\includegraphics[width=0.24\textwidth]{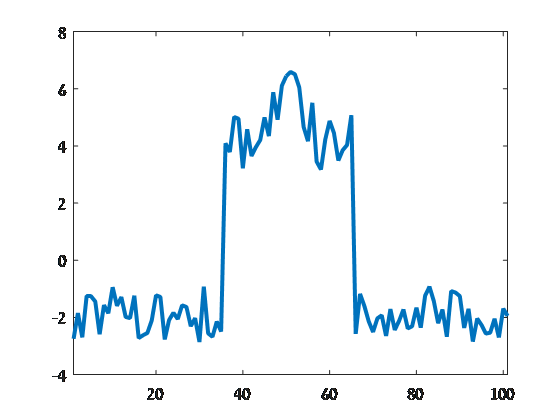} & 
	\includegraphics[width=0.24\textwidth]{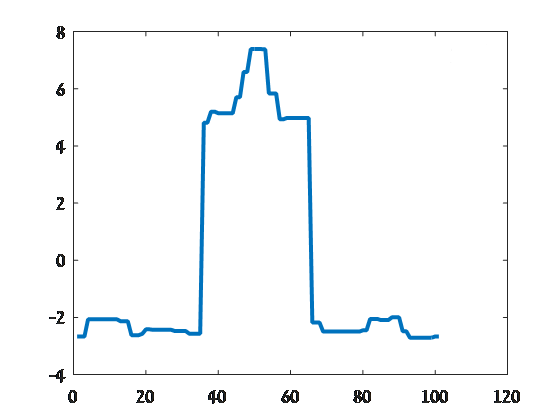} & 
	\includegraphics[width=0.24\textwidth]{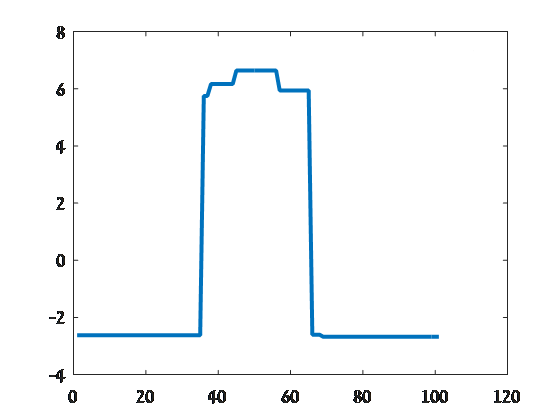} & 
	\includegraphics[width=0.22\textwidth]{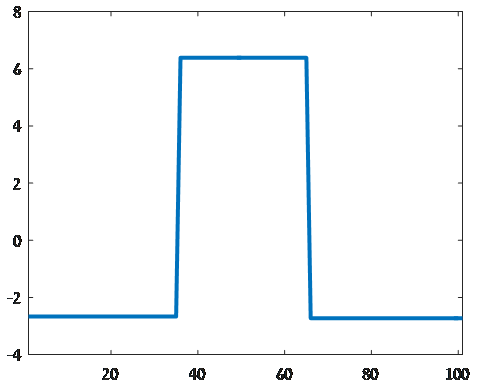} \\ 
	$u^0$ & $u^{35}$& $u^{70}$& $u^*=u^{90}$\\
	\includegraphics[width=0.22\textwidth]{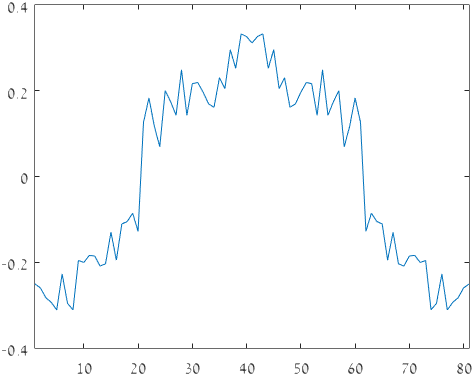}&
	\includegraphics[width=0.22\textwidth]{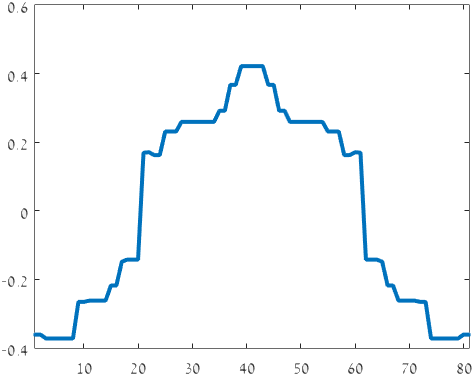}&
	\includegraphics[width=0.22\textwidth]{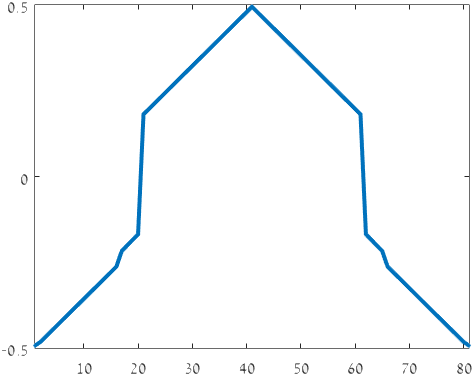}&
	\includegraphics[width=0.22\textwidth]{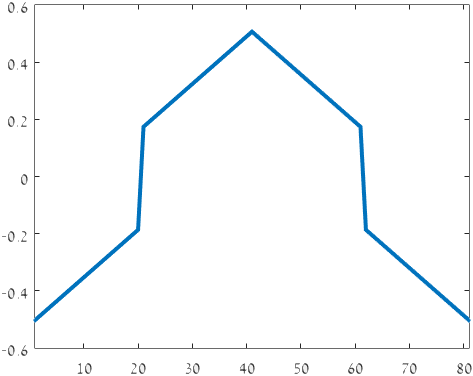}\\
		$u^0$ & $u^{4}$& $u^{8}$& $u^*=u^{12}$\\
	\end{tabular}
	\caption{Two examples of the (NG) flow. Top row $J=TV$, bottom row $J=TGV$ of order 2 (\cite{bredies_tgv_2010}).
	Taken from \cite{nossek2018flows}.}
	\label{fig:NG_flow}
\end{figure}

\begin{figure}[htb]
	\centering
\includegraphics[width=1\textwidth,]{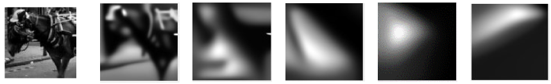}
	\caption{Nonlinear power method evolution (BHPG) for a denoising neural network FFDNet (\cite{zhang2018ffdnet}). Converged eigenfunction ($\lambda=1$), right, is a highly stable structure for the network. Taken from \cite{bungert2020powermethod}. }
	\label{fig:ffdnet}
\end{figure}

\begin{figure}[htb]
	\centering
	\begin{tabular}{cc}
\includegraphics[width=0.40\textwidth]{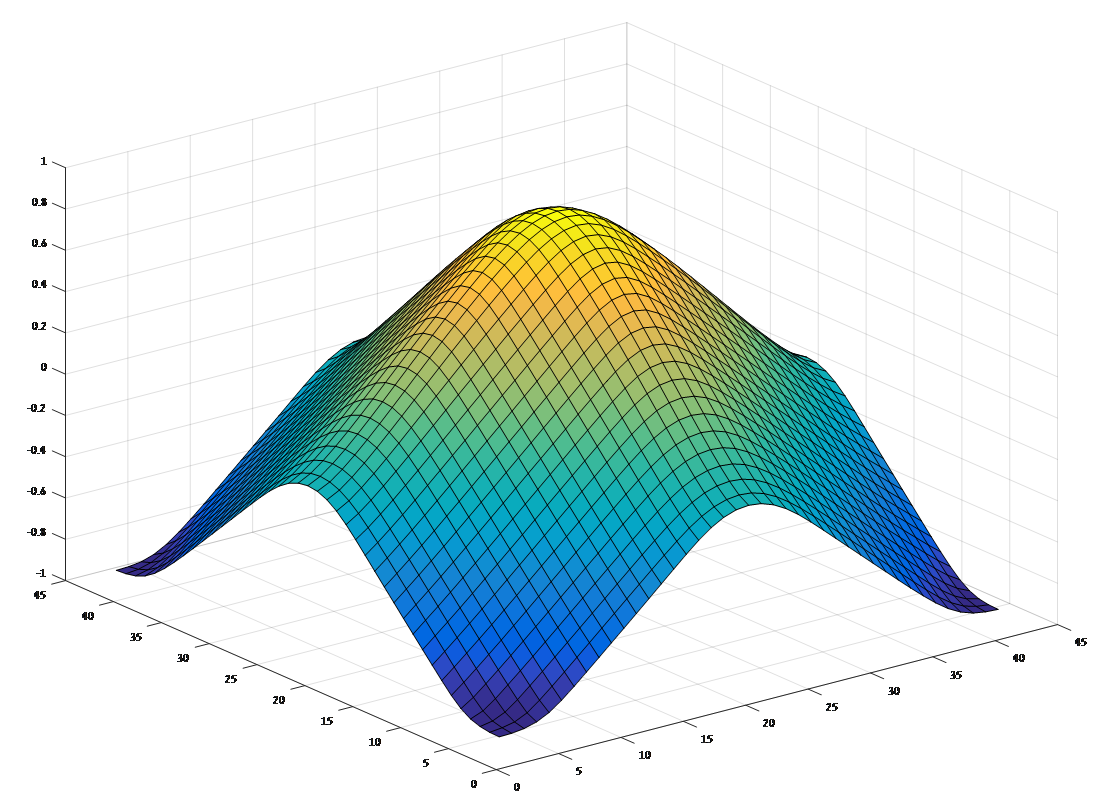}&
\includegraphics[width=0.40\textwidth]{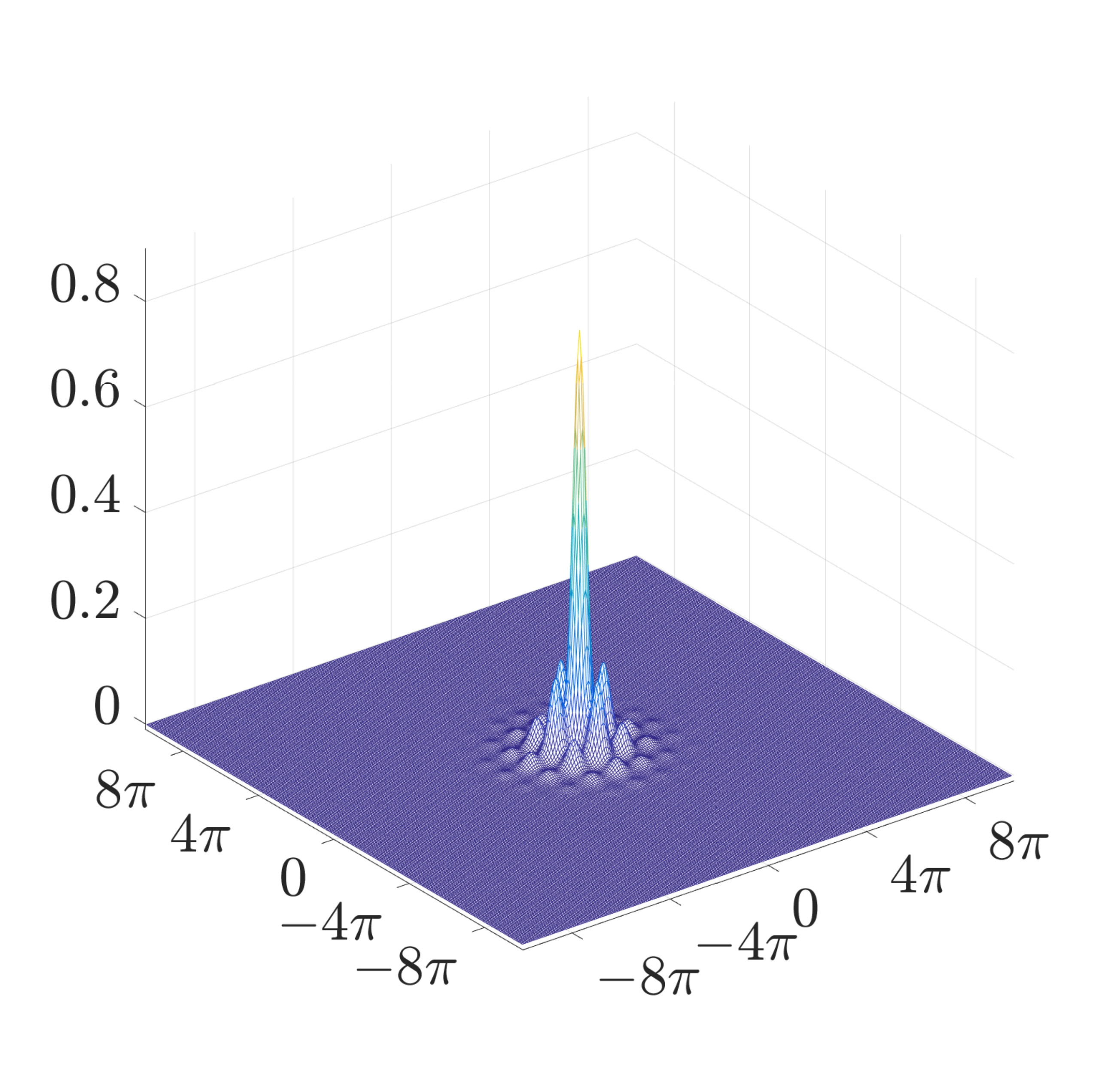}
	\end{tabular}
	\caption{EF induced by TGV (left, (NG) flow) and EF of the 2D Nonlinear Schroedinger equation \eqref{eq:nls} (right, (CG) flow). Taken from \cite{nossek2018flows} and \cite{cohen2018energy}.}
	\label{fig:EF_3dview}
\end{figure} 

%
\begin{figure}[htb]
\begin{center}
\begin{minipage}[t]{0.39\textwidth}
\begin{center}
\includegraphics[width=1\textwidth]{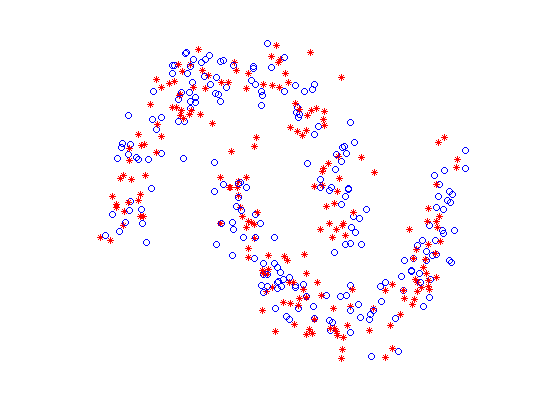}\vspace{-0.2cm}\\
Initialization
\end{center}
\end{minipage}
\begin{minipage}[t]{0.39\textwidth}
\begin{center}
\includegraphics[width=1\textwidth]{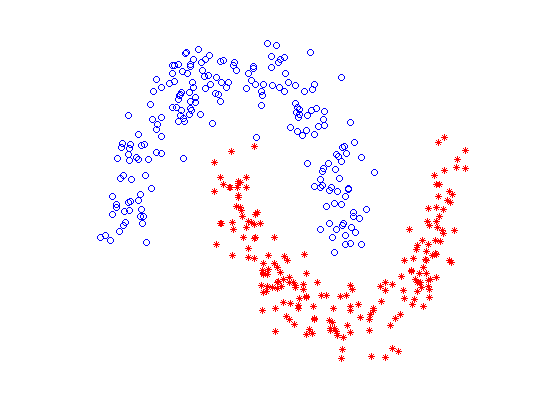}\vspace{-0.2cm}\\
Converged state
\end{center}
\end{minipage}
\caption{Results of the flow for TV defined on graphs based on point cloud distances.  The processes converges to natural clustering of the data. Taken from \cite{aujol2018theoretical}.}
\label{fig:2moons}
\end{center}
\end{figure}
 

\begin{figure}[htb]
	\centering
	\begin{tabular}{cccc}
	\includegraphics[width=0.24\textwidth]{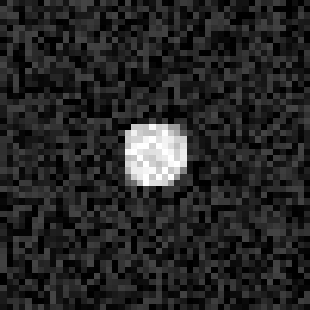} & 
	\includegraphics[width=0.24\textwidth]{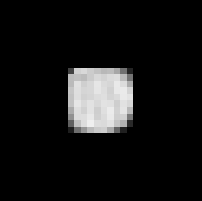} & 
	\includegraphics[width=0.24\textwidth]{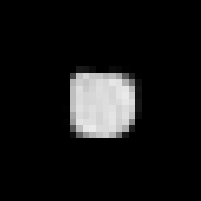} &
	\includegraphics[width=0.24\textwidth]{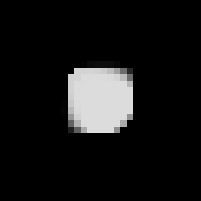} \\
	EF + noise & BM3D, PSNR=24.66$_{dB}$ & EPLL, PSNR=24.62$_{dB}$ & SpecTV, PSNR=28.12$_{dB}$
	\end{tabular}
	\caption{An eigenfunction obtained by (NG) algorithm for TV. These structure are highly stable in denoising and most suitable for the regularizer (here TV). Here it is shown that for additive white Gaussian noise, Spectral TV \cite{Gilboa_spectv_SIAM_2014} recovers well the signal, compared to well-designed classical denoisers BM3D (\cite{dabov2007image}) and EPLL (\cite{zoran2011learning}). Taken from \cite{nossek2018flows}.}
	\label{fig:EF_denoise}
\end{figure}

\begin{figure}[htb]
	\centering
	\begin{tabular}{cc}
 \includegraphics[width=0.65\textwidth]{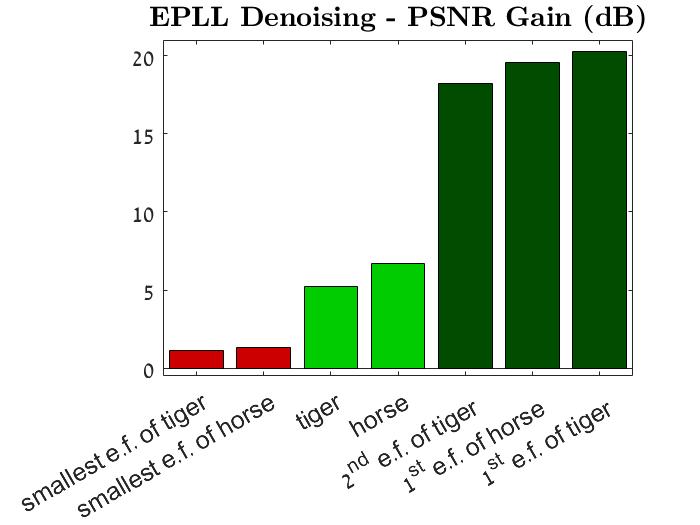}
	\end{tabular}
	\caption{Nonlinear power method for EPLL denoiser. PSNR gain: eigenfunctions vs. natural images, $var_{noise}=\frac{1}{5}var_{img}$. Taken from \cite{hait2019numeric}. }
	\label{fig:psnr_gain}
\end{figure} 

\subsection{Numerical examples}
We show some numerical examples of the algorithms presented above. In Fig. \ref{fig:NG_flow} some instances along the iteration process of (NG) are shown for the TV and TGV regularizers. At convergence we get structures which are known in the literature to be eigenfunctions induced by these functionals. 
In Fig. \ref{fig:ffdnet} we show example of the nonlinear power method (BHPG) applied to FFDNet \cite{zhang2018ffdnet}, a popular deep neural-network denoiser. We reach an eigenfunction which turns out to be a very good candidate for denoising (reaches PSNR of 44dB, compared to the horse image in the initial condition, which reached only PSNR=30dB).
In Fig. \ref{fig:EF_3dview} two examples of (NG) and (CG) flows are shown. 
Eigenfunctions on graphs are very useful for segmentation, when using graph (or nonlocal) TV for $J$, it is seen in Fig. \ref{fig:2moons} how (AGP) flow solves well the two-moon problem. Starting with a noisy initial condition (blue and red represent positive and negative values), the algorithm converges to an eigenfunction which approximate well the Cheeger cut problem.
In Figs. \ref{fig:EF_denoise} and \ref{fig:psnr_gain} we show the resilience of eigenfunctions against noise, esp. when denoised by the matching regularizer $J$ or operator $T$.
In Fig. \ref{fig:EF_denoise} an eigenfuntion of TV was denoised using 3 classical algorithms. Spectral TV (\cite{Gilboa_spectv_SIAM_2014}), which is based on the TV regularizer, is most fit to denoise such functions.
In Fig. \ref{fig:psnr_gain} we see a similar trend for EPLL denoiser. Here we have the most stable and unstable eigenfunctions (depending on their eigenvalues) and results of natural images, which are in between, with respect to denoising results. This gives insight on the priors of the denoiser, with respect to the expected spatial structures. Also adversarial examples can be obtained.

\section{Discussion and Open Problems}
In this chapter several methods for solving nonlinear eigenvalue problems are presented. Such problems appear in wide and diverse fields of signal and image processing, classification and learning and nonlinear physics. It is shown how some fundamental concepts of linear eigenvalue problems carry out to the nonlinear case. Specifically, the generalized Rayleigh quotient is a key notion, where eigenfunctions serve as its critical points.
A common theme of the presented algorithms is the use of an (often long) iterative process to compute a single eigenfunction. The process can sometimes be understood as a discrete realization of a continuous nonlinear PDE. These nonlinear flows may emerge as gradient descent of a certain energy. However, this energy is always non-convex and has many local minima (each of them is an eigenfunction). Naturally, this implies that the selection of the initial condition is critical to the computation. This is actually true for all iterative processes presented here, even if they are not directly based on a non-convex energy.
We would like to highlight several challenges this emerging field is still facing with.

We list below the main intriguing issues and open problems:
\begin{enumerate}
    \item \textbf{Initial condition.} What are the effects of the initial condition to the computation process? Can a link be formulated between the initial condition and the obtained eigenfunction? Is it related to a decomposition of the initial condition into eigenfunctions, in an analogue manner to the linear case? Are there special characteristics to the flow when random noise serves as initial condition? Is noise a good choice and in what sense?
    \item \textbf{Mapping the eigenfunction landscape of a nonlinear operator.} Can one characterize analytically eigenfunctions for a broad family of operators. This was successfully performed for TV (mainly in 2D). For more complex operators and complicated domains or graphs, this is still an open problem. For a given operator, how to design numerically algorithms which span well its eigenfunctions? We have shown that eigenfunctions of large and small eigenvalues can be computed, however reaching middle-range eigenvalues is highly non-trivial without prohibitively large computational efforts (passing through all eigenvalues in ascending/descending order).
    \item \textbf{Spectral decomposition.} Can a general theory be developed related to the decomposition of a signal into nonlinear eigenfunctions? For the case of one-homogeneous functionals, it was shown how gradient descent flows can be used for decomposition  (see \cite{Gilboa_spectv_SIAM_2014,burger2016spectral,bungert2019nonlinear}). A similar phenomenon was observed for the p-Laplacian case in \cite{cohen2020introducing}. Can this be extended to gradient descent of general convex functionals? Can these flows be used to generate multiple eigenfunctions in a much more efficient manner?
    \item \textbf{Convergence rates.} Until now the algorithms presented here did not deal with convergence rates. They are inherently quite slow, sometimes hundreds or even thousands of iterations are needed in order to numerically converge. A first analysis of the convergence rate of nonlinear power-methods for one-homogeneous functionals is in \cite{bungert2020powermethod}. This area surely requires additional focus.
    \item \textbf{Correspondence to the linear case.} It was shown that the extended definition of the Rayleigh quotient generalizes very well in the nonlinear setting. Are there additional properties related to eigenvalue analysis that can be generalized? For instance, for the power-method we know in the linear case that the method converges to the eigenfunction with the largest eigenvalue (which is part of the initial condition). We see a similar trend in the nonlinear case, where large eigenvalues are reached. Can this be formalized?
    \item \textbf{Neural networks as operators.} Last but not least, can neural networks benefit from this research field? We have shown in \cite{bungert2020powermethod} that one can treat an entire neural network (intended for denoising) as a single complex nonlinear operator and find some of its eigenfunctions. They represent highly stable and unstable modes (depending on the eigenvalue). Can additional insights be gained by analyzing eigenfunctions of deep neural networks? How can eigenfunctions be defined for classification networks? (where the input and output dimensions are very different). One direction is to develop singular value decomposition into a nonlinear setting, following the earlier work of \cite{benning2013ground}. One can also analyze eigenfunctions between layers in the net, the effect of gradient descent (or its stoachastic version) on eigenfunctions and more.
    For variational networks, the authors of \cite{effland2020variational} and \cite{kobler2020total} have shown interesting insights on the learned regularizers can be gained.
\end{enumerate}

\paragraph{Acknowledgements}
This work was supported by the European Union’s Horizon 2020 research and innovation programme under the Marie Skłodowska-Curie grant agreement No 777826, 
by the Israel Science Foundation (Grant No. 534/19) and by the Ollendorff Minerva Center.


\bibliographystyle{plain}  
\bibliography{Gilboa_refs} 

\begin{thebibliography}{10}

\bibitem{aujol2018theoretical}
Jean-Franois Aujol, Guy Gilboa, and Nicolas Papadakis.
\newblock Theoretical analysis of flows estimating eigenfunctions of
  one-homogeneous functionals.
\newblock {\em SIAM Journal on Imaging Sciences}, 11(2):1416--1440, 2018.

\bibitem{bellettini2002total}
G.~Bellettini, V.~Caselles, and M.~Novaga.
\newblock The total variation flow in $\mathbb{R}^n$.
\newblock {\em Journal of Differential Equations}, 184(2):475--525, 2002.

\bibitem{benning2013ground}
Martin Benning and Martin Burger.
\newblock Ground states and singular vectors of convex variational
  regularization methods.
\newblock {\em Methods and Applications of Analysis}, 20(4):295--334, 2013.

\bibitem{bozorgnia2016convergence}
Farid Bozorgnia.
\newblock Convergence of inverse power method for first eigenvalue of p-laplace
  operator.
\newblock {\em Numerical Functional Analysis and Optimization},
  37(11):1378--1384, 2016.

\bibitem{bozorgnia2019approximation}
Farid Bozorgnia.
\newblock Approximation of the second eigenvalue of the $ p $-laplace operator
  in symmetric domains.
\newblock {\em arXiv preprint arXiv:1907.13390}, 2019.

\bibitem{bredies_tgv_2010}
K.~Bredies, K.~Kunisch, and T.~Pock.
\newblock Total generalized variation.
\newblock {\em SIAM Journal on Imaging Sciences}, 3(3):492--526, 2010.

\bibitem{Brezis2}
H.~Brezis.
\newblock {\em Opérateurs maximaux monotones et semi-groupes de contractions
  dans les espaces de Hilbert}.
\newblock Norht Holland, 1973.

\bibitem{bungert2019nonlinear}
Leon Bungert, Martin Burger, Antonin Chambolle, and Matteo Novaga.
\newblock Nonlinear spectral decompositions by gradient flows of
  one-homogeneous functionals.
\newblock {\em To appear in Analysis \& PDE}, 2019.

\bibitem{bungert2019computing}
Leon Bungert, Martin Burger, and Daniel Tenbrinck.
\newblock Computing nonlinear eigenfunctions via gradient flow extinction.
\newblock In {\em International Conference on Scale Space and Variational
  Methods in Computer Vision}, pages 291--302. Springer, 2019.

\bibitem{bungert2020powermethod}
Leon Bungert, Ester Hait-Fraenkel, Nicolas Papadakis, and Guy Gilboa.
\newblock Nonlinear power method for computing eigenvectors of proximal
  operators and neural networks.
\newblock {\em arXiv preprint arXiv:2003.04595}, 2020.

\bibitem{Burger16}
M.~Burger, G.~Gilboa, M.~Moeller, L.~Eckardt, and D.~Cremers.
\newblock Spectral decompositions using one-homogeneous functionals.
\newblock {\em SIAM Journal on Imaging Sciences}, 9(3):1374--1408, 2016.

\bibitem{burger2016spectral}
Martin Burger, Guy Gilboa, Michael Moeller, Lina Eckardt, and Daniel Cremers.
\newblock Spectral decompositions using one-homogeneous functionals.
\newblock {\em SIAM Journal on Imaging Sciences}, 9(3):1374--1408, 2016.

\bibitem{cohen2018energy}
Ido Cohen and Guy Gilboa.
\newblock Energy dissipating flows for solving nonlinear eigenpair problems.
\newblock {\em Journal of Computational Physics}, 375:1138--1158, 2018.

\bibitem{cohen2020introducing}
Ido Cohen and Guy Gilboa.
\newblock Introducing the p-laplacian spectra.
\newblock {\em Signal Processing}, 167:107281, 2020.

\bibitem{dabov2007image}
Kostadin Dabov, Alessandro Foi, Vladimir Katkovnik, and Karen Egiazarian.
\newblock Image denoising by sparse 3-d transform-domain collaborative
  filtering.
\newblock {\em IEEE Transactions on image processing}, 16(8):2080--2095, 2007.

\bibitem{effland2020variational}
Alexander Effland, Erich Kobler, Karl Kunisch, and Thomas Pock.
\newblock Variational networks: An optimal control approach to early stopping
  variational methods for image restoration.
\newblock {\em Journal of Mathematical Imaging and Vision}, pages 1--21, 2020.

\bibitem{feld2019rayleigh}
Tal Feld, Jean-Fran{\c{c}}ois Aujol, Guy Gilboa, and Nicolas Papadakis.
\newblock Rayleigh quotient minimization for absolutely one-homogeneous
  functionals.
\newblock {\em Inverse Problems}, 35(6):064003, 2019.

\bibitem{gautier2020computing}
Antoine Gautier, Matthias Hein, and Francesco Tudisco.
\newblock Computing the norm of nonnegative matrices and the log-sobolev
  constant of markov chains.
\newblock {\em arXiv preprint arXiv:2002.02447}, 2020.

\bibitem{gautier2019perron}
Antoine Gautier, Francesco Tudisco, and Matthias Hein.
\newblock The perron--frobenius theorem for multihomogeneous mappings.
\newblock {\em SIAM Journal on Matrix Analysis and Applications},
  40(3):1179--1205, 2019.

\bibitem{Gilboa_spectv_SIAM_2014}
G.~Gilboa.
\newblock A total variation spectral framework for scale and texture analysis.
\newblock {\em SIAM Journal on Imaging Sciences}, 7(4):1937--1961, 2014.

\bibitem{gilboa2013spectral:25}
Guy Gilboa.
\newblock A spectral approach to total variation.
\newblock In {\em International Conference on Scale Space and Variational
  Methods in Computer Vision}, pages 36--47. Springer, 2013.

\bibitem{gilboa2018book}
Guy Gilboa.
\newblock {\em Nonlinear Eigenproblems in Image Processing and Computer
  Vision}.
\newblock Springer, 2018.

\bibitem{hait2019numeric}
Ester Hait-Fraenkel and Guy Gilboa.
\newblock Numeric solutions of eigenvalue problems for generic nonlinear
  operators.
\newblock {\em arXiv preprint arXiv:1909.12775}, 2019.

\bibitem{hein2010inverse}
Matthias Hein and Thomas B{\"u}hler.
\newblock An inverse power method for nonlinear eigenproblems with applications
  in 1-spectral clustering and sparse pca.
\newblock In {\em Advances in Neural Information Processing Systems}, pages
  847--855, 2010.

\bibitem{kobler2020total}
Erich Kobler, Alexander Effland, Karl Kunisch, and Thomas Pock.
\newblock Total deep variation: A stable regularizer for inverse problems.
\newblock {\em arXiv preprint arXiv:2006.08789}, 2020.

\bibitem{Meyer[1]}
Y.~Meyer.
\newblock Oscillating patterns in image processing and in some nonlinear
  evolution equations, March 2001.
\newblock The 15th Dean Jacquelines B. Lewis Memorial Lectures.

\bibitem{nossek2018flows}
Raz~Z Nossek and Guy Gilboa.
\newblock Flows generating nonlinear eigenfunctions.
\newblock {\em Journal of Scientific Computing}, 75(2):859--888, 2018.

\bibitem{rof92}
L.~Rudin, S.~Osher, and E.~Fatemi.
\newblock Nonlinear total variation based noise removal algorithms.
\newblock {\em Physica D}, 60:259--268, 1992.

\bibitem{BressonSzlam2010Cheeger}
A.D. Szlam and X.~Bresson.
\newblock Total variation and {C}heeger cuts.
\newblock In {\em International Conference on Machine Learning (ICML'10)},
  pages 1039--1046, 2010.

\bibitem{Apidopoulos}
Apidopoulos Vassilis, Aujol Jean-Fran{\c{c}}ois, and Charles Dossal.
\newblock The differential inclusion modeling fista algorithm and optimality of
  convergence rate in the case b $\backslash$leq3.
\newblock {\em SIAM Journal on Optimization}, 28(1):551--574, 2018.

\bibitem{zabusky1965interaction}
Norman~J Zabusky and Martin~D Kruskal.
\newblock Interaction of" solitons" in a collisionless plasma and the
  recurrence of initial states.
\newblock {\em Physical review letters}, 15(6):240, 1965.

\bibitem{zhang2018ffdnet}
Kai Zhang, Wangmeng Zuo, and Lei Zhang.
\newblock {FFDNet}: Toward a fast and flexible solution for {CNN}-based image
  denoising.
\newblock {\em IEEE Transactions on Image Processing}, 27(9):4608--4622, 2018.

\bibitem{zoran2011learning}
Daniel Zoran and Yair Weiss.
\newblock From learning models of natural image patches to whole image
  restoration.
\newblock In {\em Int. Conf. on Computer Vision}, pages 479--486. IEEE, 2011.

\end{thebibliography}


\end{document}